\documentclass{gtart_h}

\def\ifplaintex{\expandafter\ifx\csname documentclass\endcsname\relax}

\def\gtp{{\mathsurround=0pt\it $\cal G\mskip-2mu$eometry \&\ 
$\cal T\!\!$opology $\cal P\!$ublications}}  

\def\recd{{\small Received:\qua\receiveddate\ifx\reviseddate\relax
\else\qquad Revised:\qua\reviseddate\fi\par}} 

\def\lognumber#1{\def\thelognumber{#1}}
\def\volumenumber#1{\def\thevolumenumber{#1}}
\def\volumeyear#1{\def\thevolumeyear{#1}}
\def\papernumber#1{\def\thepapernumber{#1}}
\def\pagenumbers#1#2{\def\startpage{#1}\def\finishpage{#2}}
\def\published#1{\def\publishdate{#1}}

\def\received#1{\def\receiveddate{#1}}

\def\accepted#1{\def\accepteddate{#1}}

\def\asciiaddress#1{\def\theasciiaddress{#1}}
\def\asciiemail#1{\def\theasciiemail{#1}}

\long\def\asciiabstract#1{\long\def\theasciiabstract{#1}}
\def\asciikeywords#1{\def\theasciikeywords{#1}}

\let\\\par\let\thelognumber\relax\let\thevolumenumber\relax
\let\thepapernumber\relax\let\thevolumeyear\relax\let\startpage\relax
\let\finishpage\relax\let\publishdate\relax\let\receiveddate\relax
\let\reviseddate\relax\let\accepteddate\relax\let\theasciititle\relax
\let\theasciiauthors\relax\let\theasciiaddress\relax
\let\theasciiabstract\relax\let\theasciikeywords\relax

\let\theasciiemail\relax

\ifplaintex
\font\logobig=cmssbx10 scaled 3836
\font\logomed=cmssbx10 scaled 2557
\else
\font\logobig=cmssbx10 scaled 4200
\font\logomed=cmssbx10 scaled 2800
\fi

\long\def\makeagttitle{   
\count0=\startpage
\agt\hfill      
\hbox to 45truept{\vbox to 0pt{\vglue -13truept{\logomed A\kern -.37em{\logobig 
T}\kern -.38em G}\vss}\hss}
\break
{\small Volume \thevolumenumber\ (\thevolumeyear)
\startpage--\finishpage\nl
Published: \publishdate}

\vglue .25truein

{\parskip=0pt\leftskip 0pt plus
1fil\def\\{\par\smallskip}{\Large\bf\thetitle}\par\medskip} \vglue
0.05truein

{\parskip=0pt\leftskip 0pt plus 1fil\def\\{\par}{\sc\theauthors}
\par\medskip}%
 
\vglue 0.03truein 

{\small\leftskip 25truept\rightskip 25truept{\bf Abstract}\stdspace\theabstract

{\bf AMS Classification}\stdspace\theprimaryclass
\ifx\thesecondaryclass\relax\else; \thesecondaryclass\fi\par
{\bf Keywords}\stdspace \thekeywords\par}\vglue 7truept

}   

\ifplaintex
\font\phead=cmsl9 scaled 950
\font\pnum=cmbx10 scaled 913
\font\pfoot=cmsl9 scaled 950
\headline{\vbox to 0pt{\vskip -4.5mm\line{\small\phead\ifnum
\count0=\startpage ISSN 1472-2739 (on-line) 1472-2747 (printed)
\hfill {\pnum\folio}\else\ifodd\count0\def\\{ }%
\ifx\theshorttitle\relax\thetitle\else\theshorttitle\fi\hfill{\pnum\folio}
\else\def\\{ and }{\pnum\folio}\hfill\ifx\theshortauthors\relax\theauthors
\else\theshortauthors\fi\fi\fi}\vss}}
\footline{\vbox to 0pt{\vglue 0mm\line{\small\pfoot\ifnum\count0=\startpage
\copyright\ \gtp\hfill\else
\agt, Volume \thevolumenumber\ (\thevolumeyear)\hfill\fi}\vss}}
\else
\font=cmsl9 scaled 1050
\font\lnum=cmbx10 
\font=cmsl9 scaled 1050
\makeatletter
\def\@oddhead{{\small\lhead\ifnum\count0=\startpage ISSN 1472-2739 
(on-line) 1472-2747 (printed)\hfill {\lnum\number\count0}\else\ifodd\count0
\def\\{ }\ifx\theshorttitle\relax \thetitle \else\theshorttitle\fi\hfill
{\lnum\number\count0}\else\def\\{ and }{\lnum\number\count0}
\hfill\ifx\theshortauthors\relax 
\theauthors\else\theshortauthors\fi\fi\fi}}\def\@evenhead{\@oddhead}
\def\@oddfoot{\small\lfoot\ifnum\count0=\startpage\copyright\ \gtp\hfill\else
\agt, Volume \thevolumenumber\ (\thevolumeyear)\hfill\fi}
\def\@evenfoot{\@oddfoot}
\makeatother
\fi
\let\maketitlepage\makeagttitle
\let\makeshorttitle\maketitlepage
\let\maketitle\maketitlepage

\newwrite\gtoutfile
\long\gdef\makeheadfile{  
{\def\\{, }\def\s{ }
\immediate\openout\gtoutfile head.xxx
\immediate\write\gtoutfile{Proxy-for: \ifx\theasciiauthors\relax
\theauthors\else\theasciiauthors\fi\s<\ifx\theasciiemail\relax\theemail\else\theasciiemail\fi>}
\immediate\write\gtoutfile{\noexpand\\}
\immediate\write\gtoutfile{Authors: \ifx\theasciiauthors\relax
\theauthors\else\theasciiauthors\fi}
{\def\\{ }\immediate\write\gtoutfile{Title: \ifx\theasciititle\relax
\thetitle\else\theasciititle\fi}}
\immediate\write\gtoutfile{Subj-class: GT or SG, GR etc}
\immediate\write\gtoutfile{MSC-class: \theprimaryclass\ifx\thesecondaryclass\relax\else, \thesecondaryclass\fi}
\immediate\write\gtoutfile{Journal-ref: Algebr. Geom. Topol. \thevolumenumber\s
(\thevolumeyear) \startpage-\finishpage}
\immediate\write\gtoutfile{Comments: Published by Algebraic and
Geometric Topology at}
\immediate\write\gtoutfile{\s\s\s  http://www.maths.warwick.ac.uk/agt/AGTVol\thevolumenumber/agt-\thevolumenumber-\thepapernumber.abs.html}
\immediate\write\gtoutfile{\noexpand\\}
\immediate\write\gtoutfile{}
\ifx\theasciiabstract\relax
\immediate\write\gtoutfile{\theabstract}\else
\immediate\write\gtoutfile{\theasciiabstract}\fi
\immediate\write\gtoutfile{}
\immediate\write\gtoutfile{\noexpand\\}
\immediate\write\gtoutfile{}
\immediate\closeout\gtoutfile}}  

\def\maketitlepage{\makeagttitle\makeheadfile}
\let\makeshorttitle\maketitlepage
\let\maketitle\maketitlepage

\lognumber{24}
\volumenumber{4}
\volumeyear{2004}
\papernumber{24}
\published{11 July 2004}
\pagenumbers{521}{542}
\received{10 May 2004}
\accepted{29 June 2004}

\usepackage{amsmath,verbatim,amssymb,amscd}

\theoremstyle{plain}
\newtheorem{theorem}{Theorem}[section]
\newtheorem{corollary}[theorem]{Corollary}
\newtheorem{lemma}[theorem]{Lemma}
\newtheorem{proposition}[theorem]{Proposition}

\theoremstyle{definition}
\newtheorem{definition}[theorem]{Definition}

\theoremstyle{remark}

\newcommand{\C}{\mathbb{C}}
\newcommand{\bd}{\partial}

\def\co{\colon\thinspace}

\newcommand{\mc}[1]{\mathcal{#1}}

\newcommand{\1}{1}
\newcommand{\V}{\mc{V}}

\newcommand{\7}[6]{\left|\begin{array}{ccc} #1 & #2 & #3 \\ #4 & #5 & #6\end{array}\right|}
\newcommand{\Hom}{\text{Hom}}

\begin{document}
\title[Triangulations of $3$--dimensional pseudomanifolds]{Triangulations of $3$--dimensional pseudomanifolds\\with an application to state-sum invariants   }                    
\authors{Markus Banagl\\Greg Friedman}                  
\address{Mathematisches Institut, 
Universit\"at Heidelberg\\D-69120 Heidelberg, Germany}\secondaddress{Department of Mathematics, Yale University\\New Haven, CT 06520, USA}                  
\asciiaddress{Mathematisches Institut, 
Universitat Heidelberg\\D-69120 Heidelberg, Germany\\and\\Department of Mathematics, Yale University\\New Haven, CT 06520, USA}                  
\asciiemail{banagl@mathi.uni-heidelberg.de, friedman@math.yale.edu}
\gtemail{\mailto{banagl@mathi.uni-heidelberg.de}{\qua\rm and\qua}\mailto{friedman@math.yale.edu}}

\begin{abstract}   
We demonstrate the triangulability of compact $3$--dimensional topological pseudomanifolds and study the properties of such triangulations, including the Hauptvermutung and relations by Alexander star moves and Pachner bistellar moves. We also provide an application to state-sum invariants of $3$--dimensional topological pseudomanifolds.
\end{abstract}
\asciiabstract{%
We demonstrate the triangulability of compact 3-dimensional
topological pseudomanifolds and study the properties of such
triangulations, including the Hauptvermutung and relations by
Alexander star moves and Pachner bistellar moves. We also provide an
application to state-sum invariants of 3-dimensional topological
pseudomanifolds}


\asciikeywords{Pseudomanifold, triangulation, Hauptvermutung, Alexander star move, bistellar move, Pachner move, state-sum invariant, Turaev-Viro invariant, quantum invariant }                    

\maketitle

\section{Introduction}

In the 1950s, Moise \cite{Mo52} showed that any compact $3$--dimensional topological manifold can be triangulated and that any two such triangulations are combinatorially equivalent. This essentially means that given two triangulations of the same manifold, either can be ``straightened out'' by an isotopy so as to be piecewise linear with respect to the other. This result can be combined with earlier work of J.W. Alexander \cite{Alex} to show that any triangulation of a given $3$--manifold can be obtained from any other by an isotopy and a series of \emph{star moves} that modify the triangulation locally in a well--defined manner. There remained, however, the technical difficulty that there are an infinite number of distinct star moves (since, for example, an edge in a triangulation can be a face of an arbitrary number of tetrahedra), and this makes proofs of invariance of topological properties based upon invariance under star moves difficult: there are an infinite number of cases to test. The situation was much simplified in the 1970s by Pachner \cite{Pa}, who showed that two triangulations equivalent by star moves are also equivalent by \emph{bistellar} moves, and of these there are a finite number. So, every compact $3$--manifold is triangulable, and any two triangulations can be related by an isotopy and a finite number of combinatorial moves. 

In this paper, we seek to exptend these properties to compact $3$--dimensional topological pseudomanifolds. Pseudomanifolds are not quite manifolds, but their singularities are relatively well--behaved. Such objects arise, for example, in algebraic geometry, and they form the  principal category of  geometric object on which intersection homology is defined (see \cite{GM1, GM2}); in fact, intersection homology was created for the purpose of  extending Poincar\'e duality to pseudomanifolds.  

We demonstrate that all compact $3$--dimensional topological pseudomanifolds are triangulable and that they possess certain natural $0$-- and $1$--skeletons which will be triangulated as $0$-- and $1$--dimensional subcomplexes in any triangulation. Conversely, we show that any fixed triangulation of this natural $1$--skeleton can be extended to a triangulation of the entire pseudomanifold. The Hauptvermutung also holds for $3$--dimensional pseudomanifolds, and, as for manifolds, any two triangulations can be related by an isotopy and a series of Alexander star moves. For pseudomanifolds, we cannot reduce all Alexander moves to Pachner bistellar moves, but we will see that any two triangulations of a $3$--pseudomanifold that give equivalent triangulations on the natural $1$--skeleton can be related by an isotopy and a finite sequence of bistellar moves.

As an application of this result, we show how to extend the Turaev--Viro state-sum invariant for $3$--manifolds \cite{TV} to a family of invariants on $3$--dimensional pseudomanifolds indexed by triangulations of the natural $1$--skeleton. If the pseudomanifold $X$ has only point singularities, we recover as a topological pseudomanifold invariant the piecewise--linear invariant of Barrett and Westbury \cite{BW}. We also mention how these techniques can be used to obtain knot invariants.

\section{Triangulability and combinatorial properties of 3--pseudomanifolds}\label{S: geo}

Let us begin  by defining our objects of study.

\begin{definition}
Let $c(Z)$ denote the open cone on the space $Z$,
and let $c(\emptyset)$ be a point.  A \emph{topological stratified pseudomanifold $X$ of dimension $n$} is a  paracompact Hausdorff space
$X$ possessing a filtration
\begin{equation*}
X=X_n\supset X_{n-1}=X_{n-2} \supset X_{n-3}\supset \cdots \supset X_0\supset X_{-1}=\emptyset
\end{equation*}
such that 
\begin{enumerate}
\item $S_{n-k}=X_{n-k}-X_{n-k-1}$ is either empty or a topological manifold of dimension $n-k$ (these subsets are the \emph{strata} of $X$),
\item $X-X_{n-2}$ is dense in $X$,
\item for each $x\in S_{n-k}$ there exists a \emph{distinguished neighborhood} $U$ of $x$ in $X$ which is homeomorphic to the product $D^{n-k}\times cL$, where $D^{n-k}$ is an open ball neighborhood of $x$ in $S_{n-k}$ and $L$, the \emph{link} of $x$, is a compact topological stratified pseudomanifold of dimensions $k-1$. Furthermore, the homeomorphism $h\co  U\to D^{n-k}\times cL$ is compatible with the stratification so that is takes $U\cap X_{n-j}$ homeomorphically onto $D^{n-k}\times cL_{k-j-1}$.   
\end{enumerate}
The subspace $X_{n-2}$ is often referred to as the \emph{singular locus} and denoted by $\Sigma$.  Note that we do not allow  $X_{n-1}\neq X_{n-2}$, so, in particular, there is no $n-1$-dimensional stratum $S_{n-1}$.  

$X$ is called a \emph{topological pseudomanifold} if it can be endowed with the structure of a topological stratified pseudomanifold via some choice of filtration. 
\end{definition}

In case $X$ can be triangulated and is given a fixed triangulation, this definition agrees with the traditional definition of a pseudomanifold as an $n$--dimensional complex each of whose $n-1$ simplices bounds exactly two $n$--simplices. We will prove this below for the case $n=3$. Note, however, that there need not be a path from any $n$--simplex to any other that passes only through the interiors of $n$-- and $n-1$ simplices. 

It is also important to observe that the filtration of a pseudomanifold is, in general, not unique and that the stratum in which a point lies does not give full information about its local topology. For example, a sphere can be stratified in many ways. On the other hand, there are certain ``natural'' stratifications, which we shall discuss below.

For convenience, in what follows we will often simply write ``pseudomanifold'' to mean a topological stratified pseudomanifold, particularly when the particular choice of stratification is not relevant.

\begin{theorem}
Let $X$ be a compact topological stratified pseudomanifold of dimension $\leq 3$. Then $X$ is triangulable. In particular, there exists a finite simplicial complex $K$ and a homeomorphism $h\co  |K| \to X$. 
\end{theorem}

\begin{proof}
Let us assume that $X$ is connected. Otherwise, the following arguments can be applied separately to each connected component. 

The statement concerning triangulability will follow from a result of Munkres \cite{MK3} according to which any locally triangulable space of dimension $\leq 3$ is triangulable. We require a few definitions:

\begin{definition}
Let $U$ be an open subset of $X$, $J$ a finite complex, and $h$ a homeomorphism from $U$ to a subset of $|J|$ such that $|J|=\overline{h(U)}$ and $|J|-h(U)$ is the polytope of a subcomplex of $J$. Then $(U,h,J)$ is called a \emph{polyhedral neighborhood} on $X$.  
\end{definition}

\begin{definition}
A space is called \emph{locally polyhedral} if it is a separable metric space endowed with a \emph{locally polyhedral structure}. 
\end{definition}

For our purposes, it is enough to know that a space has a \emph{locally polyhedral structure} if it can be covered by a set, $A$, of polyhedral neighborhoods such that if $x$ and $y$ are two points in the neighborhood $(U,h,J)\in A$ and they lie in the  same open simplex of $U$ (meaning that $h(x)$ and $h(y)$ belong to the same open simplex in $J$) then 
\begin{multline*}
\text{max$\{$dim$(s)|(V,k, K)\in A$, $s$ is an open simplex in $K$, and $k(x)\in s\}$}\\
\text{$=$max$\{$dim$(s)|(V,k, K)\in A$, $s$ is an open simplex in $K$, and $k(y)\in s\}$}.
\end{multline*}
In other words, if two points share a simplex under some local triangulation of the covering $A$, then they must have the same \emph{index} with respect to the covering, where the index of a point (with respect to $A$) is the maximum dimension of  the open simplices in which it lies, taken over all polyhedral neighborhoods in the collection $A$. 
The dimension of a polyhedral space is the highest dimension of a simplex which occurs in any of the polyhedral neighborhoods. N.B. the index of a point is determined by the structure of the family of polyhedral neighborhoods and \emph{not} necessarily by actual local topology. 

We seek then to employ Theorem 9.8 of Munkres \cite{MK3}, which states that any locally polyhedral space of dimension $n\leq 3$ can be globally triangulated. Hence we must show that a compact topological stratified pseudomanifold $X$ of dimension $\leq 3$ is a locally polyhedral space.

It will be convenient to perform an induction, so let us begin by noting that a compact pseudomanifold of dimension $0$ must be a finite collection of points. Clearly, this has a triangulation. It follows similarly from the definitions that a compact pseudomanifold of dimension $1$ must be a disjoint union of circles, also triangulable. So now assume inductively that we have established triangulability of pseudomanifolds of all dimensions $\leq n-1\leq 2$. 

Let $B$ denote the covering of $X$ given by all distinguished neighborhoods of points of  $X$ (this is a covering since each point has a distinguished neighborhood). We will show that $B$ determines a locally polyhedral structure (although $B$ itself may only be the basis for the locally polyhedral structure; see \cite{MK3}). First, we must show that each distinguished neighborhood determines a polyhedral neighborhood. Let $U$ be some distinguished neighborhood, and recall that there is a homeomorphism  $h\co  U\to D^{n-k}\times cL$. Clearly, $\bar D^{n-k}$ is triangulable (in fact by the $n-k$ simplex, $\bar \sigma^{n-k}$). We can assume that the link $L$ has been triangulated by induction, and this induces a triangulation on the closed cone $\bar cL$: if $\mc{L}$ is a simplicial complex such that $g\co  L\to \mc{L}$ is a homeomorphism, then we define the homeomorphism $\bar cg\co \bar cL\to \bar c\mc{L}$ by $\bar cg(y,t)=(g(y),t)$. Therefore, since $\bar \sigma^{n-k}$ is obviously a locally finite complex, $|\bar \sigma^{n-k}|\times |\bar c\mc{L}|$ can be triangulated by a complex $J$  such that if $\bar \sigma$ is a closed simplex in $\bar \sigma^{n-k}$ and $\bar \tau$ is a closed simplex in $c\mc{L}$, then each $|\sigma| \times |\tau|$ is the polytope of a subcomplex of the triangulation of $|\bar \sigma^{n-k}|\times |\bar c\mc{L}|$
(see \cite[p.\ 339]{MK}). But  $|J|=|\bar \sigma^{n-k}|\times |\bar c\mc{L}|\cong \bar D^{n-k}\times \bar cL$ by a homeomorphism $f$ that clearly takes $D^{n-k}\times cL$ to a set whose closure is all of  $|J|=|\bar \sigma^{n-k}|\times |\bar c\mc{L}|$ and such that   $|J|-f( D^{n-k}\times cL)=(|\bd \bar \sigma^{n-k}|\times |L|)\cup (|\bar \sigma^{n-k}|\times |\mc{L}|)$ which is the polytope of a subcomplex of $J$. Together, this shows that $(U, fh, J)$ gives a polyhedral neighborhood of $x$.

Next, we need to check the condition on the indices of points under these coverings by distinguished neighborhoods which, a posteriori, are coverings by polyhedral neighborhoods. In order to facilitate this discussion, let us be more specific about our choices of local triangulations for the distinguished neighborhoods in the cover (if necessary, we alter the covering to obtain a new one with the same sets $U$ covering $X$ but different choices of triangulating complex and homeomorphism). Due to our limitations in dimension, there are only three possibilities for distinguished neighborhoods: 
\begin{enumerate}
\item \label{I: top nghbd} $U$ is a distinguished neighborhood of a point in the top stratum $X-X_{n-2}$. In this case, $U\cong D^n$, and we can choose $J=\bar \sigma^{n}$ and $f\co D^n\to |\sigma^n|$ to be any of the standard homeomorphisms  induced by a homeomorphism $\bar D^n\to |\bar \sigma^n|$. Note that by the stratum compatibility conditions for distinguished neighborhoods, $U\cap X_{n-2}=\emptyset$. 

\item \label{I: 1 nghbd} $U$ is a distinguished neighborhood of a point in $X_1$ (in which case $X$ must have dimension $3$). Then $U\cong D^1\times cL$. In this case, we note that it is possible to choose the triangulation of the polyhedron $|\bar\sigma^1|\times |\bar c\mc{L}|$ so that $|\bar \sigma^1|\times |*|$, where $*$ is the cone point of $\bar c\mc{L}$, is triangulated by a subcomplex consisting of a single closed $1$--simplex (see \cite[p.\ 104--105]{MK}). If we form the polyhedral neighborhood using a triangulation of this type, then $D^1\times *$ will be mapped homeomorphically onto an open $1$--simplex of $J$. Note that by the stratum compatibility conditions for distinguished neighborhoods, $U\cap X_{0}=\emptyset$. 

\item \label{I: 0 nghbd} $U$ is a distinguished neighborhood of a point $x$ in $X_0$. In this case $U\cong cL$, where $L$ is a pseudomanifold of dimension $\leq 2$ and $x=U\cap X_0$ will map to the cone point of the complex $|\bar c\mc{L}|$. Recall that we have obtained a triangulation $g\co |\mc{L}|\to L$ by an inductive application of Munkres's theorem. By a more careful application of Munkres's theorem, we can assume that $L_0$  is triangulated as a subcomplex. This is because \cite[Cor.\ 9.7]{MK3} provides a \emph{proper} triangulation of $L$ and a proper triangulation induces a triangulation of $L^0$, the set of points in $L$ with index $0$. But with the locally polyhedral structure of $L$ also determined by distinguished neighborhoods, this is exactly the set $L_0$: any point in $L-L_0$ has index $1$ or $2$ (as $\dim(X)=3$ or $2$) as we can see from item \ref{I: top nghbd} applied to $L$, while points in $L_0$ must have index $0$ by the inductive application of this item, which imposes cone--like triangulations on distinguished neighborhoods of points in $L_0$.  Note that any point in $L_0$ cannot lie in any distinguished neighborhood of any point in $L-L_0$ by the strata compatibility conditions for distinguished neighborhoods. 

Thus, we are free to assume that the image of $L_0$ is a zero dimensional subcomplex of $\mc{L}$. We then take the obvious triangulation of $\bar c\mc{L}$ as the complex underlying the triangulation of $cL$. With this choice of triangulation, the stratum compatibility conditions for distinguished neighborhoods ensure that $U\cap (X_1-X_0)$  is triangulated by a disjoint union of open $1$--simplices. 
\end{enumerate}

With these choices of triangulations for polyhedral neighborhoods, we can now demonstrate that this system of neighborhoods satisfies the needed index condition on $X$. Recall that we must show that if $x$ and $y$ are two points in the neighborhood $(U,h,J)$ and they lie in the  same open simplex of $U$ then the indices of $x$ and $y$ agree in $X$, where the index $I(x)$ is defined as the maximum dimension of the open simplices in which $x$ resides in all of the local polyhedral neighborhoods of the covering:
$I(x)=$max$\{$dim$(s)|(V,k, K)\in B$, $s$ is an open simplex in $K$, and $k(x)\in s\}$. The proof comes down to checking cases using the choices of local polyhedral neighborhoods we have just established:

\begin{enumerate}

\item If $x, y\in X$ share an open $2$--simplex or $3$--simplex in any polyhedral neighborhood in $B$, then by the above choices of local triangulations, $x$ and $y$ must lie in the top stratum $X-X_{n-2}$. Therefore, $I(x)=I(y)=n$ by looking at the distinguished neighborhoods of these points.

\item If $x,y\in X$ share an open $1$--simplex, then by the above choices of local triangulations, they are either both elements of $X-X_{n-2}$, in which case the index of each is $n$ as in the last item, or they are both elements of $X_1$, in which case they each have index $1$, again by our choice of triangulations of distinguished neighborhoods.

\item It is impossible for two distinct points $x,y \in X$ to share a $0$ simplex, so there is nothing to check here.
\end{enumerate}

Lastly, in order to be able to apply Munkres's theorem and complete our proof of triangulability, we must show that $X$ is separable and metrizable (once it is metrizable, we can then endow it with any allowable metric we'd like in order to make $X$ a separable metric space). But clearly $X$ is locally metrizable since we can endow the neighborhood of any point with the metric determined by its homeomorphic image in a finite complex, which itself possesses a metric induced from some embedding in euclidean space. Then since $X$ is compact and Hausdorff, it is metrizable. Furthermore, since $X$ is metrizable, separability is equivalent to second countability, which follows for $X$ by compactness and metrizability (see \cite{MK2}). 
\end{proof}

We next establish that there are certain natural  subsets of a 3--pseudomanifold that will be the polytopes of subcomplexes of any triangulation. This requires some preliminary work.

\begin{lemma}[Classification of 2--pseudomanifolds]
Any compact  topological pseudomanifold of dimension $2$ consists of a collection of compact $2$--manifolds without boundary glued or self--glued along a finite number of points. More technically: such a pseudomanifold $X$ is the quotient space of a disjoint union of compact $2$--manifolds without boundary such that if $q$ is the quotient map, then all points $x\in X$ satisfy $\#\{q^{-1}(x)\}<\infty$ and for all but finitely many $x\in X$, $\#\{q^{-1}(x)\}= 1$.
\end{lemma}
\begin{proof}
Let $X$ be a  compact $2$--dimensional topological pseudomanifold filtered by $X\supset X_0$. Since $X$ is compact, the strata compatibility condition for distinguished neighborhoods in the definition of a pseudomanifold ensures that $X_0$ is a finite set. We also know, from the definitions, that $X-X_0$ is an open $2$--manifold without boundary. For each $x\in X_0$, consider a distinguished neighborhood $U_x$ given by $U_x\cong cL$, where $L$ is a compact $1$--dimensional pseudomanifold. It is clear from the definitions that $L$ must be a disjoint union of circles and so $U_x\cong c(\coprod S^1)$. Consider the smaller distinguished neighborhoods $V_x\cong c(\coprod S^1)\subset U_x$ given by $V_x=\{(x,t)\in U_x| t<1/2\}$, where we here identify $cZ$ with the set of points $(z,t)$ modulo the identification $(z_1,0)\sim (z_2,0)$. Clearly then $X-\cup_{x\in X_0} V_x$ is a $2$--manifold $M$ with  a collared boundary. Now reattach the  neighborhoods $V_x$.  This can be viewed as the two step process of filling in the boundary of $M$ with disks to form a compact $2$--manifold and then attaching some of these disks appropriately along the cone points. This  implies the lemma by taking the quotient.
\end{proof}

\begin{lemma}\label{L: neigh class}
All distinguished neighborhoods of points in a $3$--pseudomanifold  $X$ have one of the following homeomorphism types:
\begin{enumerate}
\item \label{I: type 1}$D^3$,
\item \label{I: type 2}$D^1\times c(\coprod S^1)$,
\item \label{I: type 3} $cL$, where $L$ is a compact $2$--dimensional pseudomanifold; in particular a $2$--pseudomanifold as classified by the previous lemma.
\end{enumerate}
\end{lemma}
\begin{proof}
This is clear from the definition of pseudomanifold and the classification of $1$ and $2$ dimensional pseudomanifolds.
\end{proof}

Note that neighborhoods of type 1 of the lemma are also of type 2, and those of type 2 are also of type 3.

\begin{proposition}\label{P: classic def}
Suppose that the compact $3$--pseudomanifold $X$ is triangulated by the complex $K$. Then $K$ is a simplicial pseudomanifold in the sense that $K$ is a union of $3$--simplices such that every $2$--simplex in $K$ is the face of exactly two $3$--simplices (note, however, that we do not require the classical condition that every $3$--simplex can be connected to every other $3$--simplex by a path of $3$--simplices connected by $2$--faces).
\end{proposition}
\begin{proof}
That $X$ is the union of $3$--simplices follows from the density of the $3$--manifold $X-X_2$ in $X$.\

Suppose now that $\bar \sigma^2$ is any closed $2$--simplex in $X$. Note that if $\bar \sigma^2$ is the boundary of exactly two $3$--simplices, then any point of the open simplex $\sigma^2$ has a neighborhood of type \eqref{I: type 1} as defined in the previous lemma. Suppose, on the other hand that $\bar \sigma^2$ is the face of either one or more than two $3$--simplices (a priori it can not be the face of zero $3$--simplices, since $X$ is the union of $3$--simplices). But some geometric thought shows that this is inconsistent with any point in $\sigma^2$ having a distinguished neighborhood of types (1), (2), or (3) of the previous lemma. 
\end{proof}

\begin{definition}
Suppose $X$ is a compact $3$--pseudomanifold. We define two natural subsets of $X$ as follows:
\begin{enumerate}
\item Let $X(1)$ be the subset of $X$ consisting of points that have no neighborhood in $X$ homeomorphic to $D^3$ (i.e.\ of type \eqref{I: type 1} of Lemma \ref{L: neigh class}).
\item Let $X(0)$ be the subset of $X(1)$ consisting of points which also have no neighborhood homeomorphic to $D^1 \times c(\coprod S^1)$ (i.e.\ points with no neighborhood of type \eqref{I: type 1} or type \eqref{I: type 2} of Lemma \ref{L: neigh class}).
\end{enumerate}
Alternatively, $X(1)$ is the set of points of $X$ which must be in $X_1$ for \emph{any} choice of stratification of the space $X$ as a pseudomanifold, while $X(0)$ is the set of points which must be in $X_0$ for \emph{any} choice of stratification. 
\end{definition}

\begin{proposition}\label{P: subgraph}
For any triangulation $h\co X\to |K|$ of the $3$--dimensional compact pseudomanifold $X$, $h(X(0))$ and $h(X(1))$ are the polytopes of respectively $0$--dimensional and $1$--dimensional subcomplexes of $K$. 
\end{proposition}
\begin{proof}
For convenience of notation, we use the homeomorphism $h$ to identify $X$ with the polytope $K$ and identify the images of $X(0)$ and $X(1)$ with subsets (not yet necessarily subcomplexes) $K(0)$ and $K(1)$ of $|K|$. The proposition will follow from the fact that any simplicial pseudomanifold can be given the structure of a topological pseudomanifold filtered by the polytopes of the simplicial skeleta $K^i$, i.e.\ by $|K^0|\subset |K^1|\subset |K|$ for a simplicial $3$--pseudomanifold (see \cite[Chapter I]{Bo}).

First, we show that $K(0)\subset |K^0|$. Suppose to the contrary that $x\in |K|$, but $x\notin |K^0|$. Then $x$ will be  a point in one of the pure strata $|K|-|K^1|$ or $|K^1|-|K^0|$. But in the first case $x$ must have a distinguished neighborhood of type \eqref{I: type 1} and in the second it must have a neighborhood of type \eqref{I: type 2} of Lemma \ref{L: neigh class}, and these cases would contradict $x\in K(0)$. Hence $K(0)\subset |K^0|$, and clearly $K(0)$ is the polytope of a $0$--dimensional subcomplex of $K$.

Now consider $K(1)$. Similarly to the arguments of the last paragraph, we must have $K(1)\subset |K^1|$, else some point of $K(1)$ would have a neighborhood of type \eqref{I: type 1}, a contradiction. Notice also that $K(1)$ is a closed subset of $|K|$ since clearly its complement, the set of points which have neighborhoods homeomorphic to $D^3$, is open. So it now suffices to demonstrate that $K(1)$ contains every open $1$--simplex that it intersects. But for any two points $x$, $y$ in a $1$--simplex $\sigma^1$, there is a homeomorphism $|K|\to |K|$ which takes $x\to y$ (indeed, there is the standard piecewise linear homeomorphism on $\sigma^1 \to \sigma^1$ which takes $x$ to $y$ and then extends linearly to St$(\sigma^1)$ and is the identity outside this star). But then clearly $x$ and $y$ have homeomorphic neighborhoods, so if $x\in K(1)$ so must be $y$. Hence, $K(1)$ must be the polytope of the closure of a union of open $0$ and $1$ simplices, which must be a $1$--dimensional complex. 
\end{proof}

Given now  that we can triangulate any compact $3$--dimensional topological pseudomanifold and that, for each, there exists a certain naturally occurring $1$--skeleton, $X(1)$, which must appear as a subcomplex of any such triangulation, we show that there always exists a triangulation of the whole pseudomanifold that extends any given triangulation of this natural $1$--skeleton.

\begin{proposition}\label{P: realize}
Given a compact topological $3$--pseudomanifold $X$ and a triangulation of its natural $1$--skeleton $X(1)$ that is compatible with the natural stratification (in other words, each point in the $X(0)$ must be subpolyhedron, i.e.\ a vertex), there exists a triangulation of $X$ that extends the given triangulation of $X(1)$.
\end{proposition}
\begin{proof}
From Proposition \ref{P: subgraph}, we know that $X(1)$ and $X(0)$ are subpolyhedra of $X$ in any triangulation, so topologically $X(1)$ must be a graph, a collection of vertices and edges (disconnected vertices are allowed). So any triangulation of $X(1)$ is by a polygonal graph $\Gamma$, and we choose one such triangulation up to isotopy in accordance with the hypothesis of the proposition.  Fixing this triangulation only up to isotopy is not a serious restriction since any isotopy on  $X(1)$ that fixes $X(0)$ can be extended to all of $X$.  This will follow from the techniques of our proof.  So, in the end, we can pick any triangulation of $X(1)$ in the isotopy class and then isotop the remainder of the triangulation of $X$  we have constructed to suit it. Now, however, it is useful to have the added flexibility of only fixing an isotopy class.

Let us choose some  arbitrary triangulation $|K|\to X$ of $X$. For simplicity, we can identify $X$ with this complex in what follows. Of course in this triangulation $X(1)$ and $X(0)$ are the polytopes of a $1$--subcomplex and a $0$--sub--subcomplex, respectively. By an isotopy that fixes the vertices of $\Gamma$ that map to $X(0)$, we can modify the triangulation $|\Gamma|\to X(1)$ to make it piecewise linear with respect to the triangulation $K$. Then there exist subdivisions of $\Gamma$ and $K$ such that this map is simplicial. Without loss of generality (since our choice of $K$ was arbitrary to begin with), let us replace $K$ with this subdivision if necessary  (without changing its name).  With these choices,  we now see that we can consider the triangulation of $X(1)$ under $K$ to be a subdivision of that by $\Gamma$. 

Our goal now will be to show that given such a global triangulation,  we may alter it to obtain another triangulation of $X$ such that $X(1)$ is again triangulated by $\Gamma$, itself.

To fix some notation, let us identify $|\Gamma|$, the geometric realization of the graph $\Gamma$, with $X(1)$ via the given triangulation, and let $|\Gamma_0|$ denote the images of the vertices of $\Gamma$ in $X(1)$. We continue to treat $X$ as the geometric realization of  a simplicial complex $K$, and then $|\Gamma|$ and $|\Gamma_0|$ are subpolyhedra  triangulated by subcomplexes of $K$. Let us consider the second barycentric subdivision $K''$ of $K$ and the relative stellar neighborhood in $K''$ of $|\Gamma|$ rel $|\Gamma_0|$; call this neighborhood $N(\Gamma, \Gamma_0)$. In other words, $N(\Gamma, \Gamma_0)$ will be a relative regular neighborhood in $X$ of $|\Gamma|$ rel $|\Gamma_0|$ (see \cite{Co}). In concrete terms, this neighborhood is the union of (closed) simplices in $K''$ that intersect $|\Gamma|-|\Gamma_0|$. 

It is not hard to see that $N(\Gamma, \Gamma_0)$ is a union of $3$--balls. More specifically, if $|\gamma_i|$  is the geometric realization of an edge of $\Gamma$ and the link of $|\gamma_i|$ in $X$ has $n_i$ components (note that this number will be well--defined since int$(\gamma_i)$ must lie in  $X(1)-X(0)$), then the relative regular neighborhood $N(\gamma_i, \bd \gamma_i)$ will consist of $n_i$ $3$--balls joined along $\gamma_i$. In fact, take $K''$ and consider the decomposition space $K''_*$ as defined in \cite{MK3}. In this decomposition, we have simply unglued along $X(1)-X(0)$, so, in particular, each vertex of $|\Gamma_0|$ in $X(0)$  corresponds to exactly one vertex in $K''_*$ and each vertex of $|\Gamma_0|$ that lies in $X(1)$ lifts to a number of points in the decomposition space corresponding to the number of connected components in its link. We denote this set of vertices in $K''_*$ by $|\Gamma_0|_*$. Meanwhile, each $|\gamma_i|$ lifts to $n_i$ $1$--dimensional subpolyhedra that are disjoint except possibly at their endpoints in $ |\Gamma_0|_*$. We choose some ordering and label these arcs $\gamma_{ij*}$.  
Note that  $K_*''$ is a $3$--pseudomanifold with only point singularities. 

Now consider $\check K=K_*''-\cup \text{int (St}(|\Gamma_0|_*))$, i.e.\ $K''_*$ with the stars in $K''$ of the vertices in $|\Gamma_0|_*$ removed. $\check K$ is a $3$--manifold with boundary. Since we have take subdivisions, each $|\gamma_{ij}|\cap \check K$ is an arc which meets the boundary of $\check K$ normally. We have also taken sufficient subdivisions that the stellar neighborhoods of the $|\gamma_{ij*}\cap \check K|$ (i.e.\ the neighborhoods consisting of the unions of the simplices that intersect them) will be disjoint regular neighborhoods of them. By standard regular neighborhood theory in a $3$--manifold, these neighborhoods are homeomorphic to $3$--balls (see \cite{HUD}); we use here that all triangulations of $3$--manifolds are compatible with their manifold structures (this fact, originally due to Moise \cite{Mo52}, also follows from a theorem of Edward Brown quoted below as Theorem \ref{T: Brown}). Let's label these balls $M_{ij*}$. In addition, the intersection of each neighborhood $M_{i,j*}$ with the boundary is a regular $2$--disk neighborhood of $|\gamma_{ij*}|\cap \bd \check K$ (see \cite{HUD}). Now putting the stars around  $|\Gamma_0|_*$ back in, we see that the stellar neighborhood $N_{ij*}$ of $|\gamma_{ij*}|$ rel $\bd|\gamma_{ij}|$ in $K''_*$ is the union of the $M_{ij*}$ with the cones on the boundary pieces $M_{ij*}\cap \bd \check K$. Thus each $N_{ij*}$ is also  homeomorphic to a $3$--ball. Now, passing back down to $K''$ from $K_*''$, the triangulation is unaffected except for regluing the $|\gamma_{ij}|$ back together to $|\gamma_i|$. Hence each $N(\gamma_i, \bd \gamma_i)$ is homeomorphic to the union of the $N_{ij*}$ joined along $|\gamma_i|$. 
The entire relative neighborhood $N(\Gamma, \Gamma_0)$  is  the union of the $N(\gamma_i, \bd \gamma_i)$, and these are joined together only at the boundary points of the $|\gamma_i|$. Note that since $\Gamma$ is a well--defined $1$--complex, no two boundary points of the same $|\gamma_i|$ are ever joined together in this way.

Our goal now is to modify the triangulations of the neighborhoods $N(\gamma_i, \bd \gamma_i)$ so that, within each, each $|\gamma_i|$ is triangulated by a single $1$--simplex. With suitable modifications to the triangulation of $X$ off  $|\Gamma|$, we will obtain our desired triangulation. Notice that, topologically, each $N(\gamma_i, \bd \gamma_i)$ is homeomorphic to a union of standard $3$--balls, $B_{ij}$, joined along their north--south axes (this follows by an application of isotopy extension from the axis of a ball to the entire ball rel boundary). Let $\Omega_i$ denote the subcomplex of $N(\gamma_i, \bd \gamma_i)$ consisting of all simplices in $K''$ in the boundaries $\bd B_{ij}$ except for those that intersect the south pole. So $\Omega_i$ is a collection of $2$--disks $D_{ij}$ joined together at a single point of each. If we then take the closed cone $\bar c\Omega_i$, we obtain a space homeomorphic to $N(\gamma_i, \bd \gamma_i)$. Furthermore, the triangulations of $\bar c\Omega_i$ and $N(\gamma_i, \bd \gamma_i)$ agree on $\cup \bd B_{ij}$ since the open cone $c\bd D_{ij}$ is exactly what we removed from $\cup \bd B_{ij}$ to form $\Omega_i$. Lastly, of course, $\bar c\Omega_i$ has exactly one $1$--simplex that runs through its interior as the north--south axis $|\gamma_i|$. 

We perform this procedure for all $i$ and so we can remove all $N(\gamma_i, \bd \gamma_i)$ from the triangulation $K''$ and replace them with the $\bar c\Omega_i$ by gluing along the matching boundaries. This creates a $\Delta$--complex in the sense of \cite{Ha}, i.e.\ a quotient space of simplices joined along their faces. It remains to see that we can turn this into a true triangulation without further modifying $|\Gamma|$ which is now triangulated as we want it. It is clear from the construction that each simplex uniquely determines a full set of distinct vertices, but it is necessary also that each set of vertices spans a distinct (possibly empty) simplex. In other words, we must ensure that no set of vertices spans multiple simplices. One way to ensure this would be to take another subdivision of our current complex, but that would ruin what we have achieved as far as triangulating $|\Gamma|$ goes. So instead, we perform a generalized stellar subdivision rel $|\Gamma|$, a subdivision built on the current set of vertices along with the barycenters of all simplices \emph{not} in $|\Gamma|$ (see \cite{MK}). Let us call the resulting ``complex'' $L$; we show that $L$ is in fact a legitimate simplicial complex and not just a $\Delta$--complex.

Clearly, there is no trouble with $0$--simplices of $L$. Nor is there any ambiguity with simplices that have a vertex either in the interior of some $|\bar c\Omega_i|$ (by which we mean the union of the interiors of the corresponding $B_{ij}$) or on the complement of the interior of $\cup|\bar c\Omega_i|$, since the restriction of $L$ to $|\bar c\Omega_i|$ and the complement of its interior must each be simplicial complexes (as subdivisions of simplicial complexes) and there are clearly no simplices in $L$ with vertices both interior and exterior to any one $|\bar c\Omega_i|$. So we need only worry about simplices whose vertices all lie on some $|\bd  \bar c\Omega_i|$ (by which we mean the union of the corresponding $|\bd B_{ij}|$).  In fact, we need only be concerned if they all lie on the same $|\bd B_{ij}|$ since any two vertices that lie only on distinct $\bd B_{ij}$ (thus excluding the north and south poles)  cannot span any $1$--simplex either within $|\bar c\Omega_i|$, by construction, or on the complement of its interior, because of the extra generalized subdivision we have made to obtain $L$ (there can be no $1$--simplices  left which touch multiple $| B_{ij}|$ except possibly at the north or south poles). 

So, let us consider $1$--simplices whose vertices all lie on the same $|\bd B_{ij}|$. 
Due to the generalized subdivision, the only possibilities for such a $1$--simplex, $\sigma$, are for it to lie in some $|\bd B_{ij}|$ or for it to be $|\gamma_i|$: in order for $\sigma$ not to be a $\gamma_i$ and still intersect the interior of some $| B_{ij}|$, it would have to include a vertex also in the interior of $| B_{ij}|$; similarly it cannot lie in the complement of the interior of $|\bar c\Omega_i|$. 
If $\sigma=\gamma_i$, its vertices are the north and south poles of the $|B_{ij}|$. By assumption that $\Gamma$ is properly triangulated by the $\gamma_i$, these vertices cannot span any other $1$--simplex in $|\Gamma|$, and since the $|\bar c\Omega_i|$ only possibly intersect in their various north and south poles, $\bd \gamma_i$ cannot span any other $1$--simplices in the restriction of $L$ to $\cup |\bar c\Omega_i|$.  On the other hand, these vertices again cannot span any $1$--simplex in the complement of $\cup |\bar c\Omega_i|$ because of the generalized subdivision we performed there to obtain $L$ since any such spanning $1$--simplex would have been subdivided at that time (of course the subdivision could not have created an extra $1$--simplex connecting the two vertices). If $\sigma\neq \gamma_i$, then $\sigma$ lies in  some $|\bd B_{ij}|$. But since $L$ restricts to a triangulation of $|\bd B_{ij}|$, which is a subcomplex of $\cup|\bar c\Omega_i|$, $\bd \sigma$ can not span any other simplices in $|\bd B_{ij}|$. Thus we have ruled out any other possible simplices for $\bd \sigma$ to span except for $\sigma$. 

Next, let us consider  $2$-- or $3$--simplices whose vertices all lie on the same $|\bd B_{ij}|$. Once again because of the generalized subdivision, any such simplex must lie entirely in $|\bd B_{ij}|$, using the observation we have already made that any pair of vertices of any such simplex must uniquely determine its spanning $1$--simplex. In this case, none of these $1$--simplices can be $\gamma_i$, since the north and south poles of $|B_{ij}|$ are not connected by any path with two edges in  $|\bd B_{ij}|$. It then follows that any such simplex is uniquely determined in $|\bd B_{ij}|$ by using that the complement of the interior of $\cup |\bar c\Omega_i|$ is properly triangulated.

Thus, we have shown that $L$ is a legitimate simplicial complex that triangulates $X$ and imposes the desired triangulation on $X(1)$. The promised ability to extend isotopies on $X(1)$ rel $X(0)$ comes from continuing to view the $|\gamma_i|$  as axes of balls and the ability to extend such isotopies rel the boundaries of the balls.
\end{proof}

Now that we have determined that any compact $3$--pseudomanifold $X$ can be triangulated and studied both the limitations and allowable possibilities for such triangulations, let us look at relations between these various triangulations.

\begin{definition}
An \emph{Alexander star move} on an open simplex $E$ in a polyhedron triangulated by a simplicial complex $T$ consists of replacing the closed star St$(E)$ (the union of the closed simplices in $T$ which contain $E$) with the cone to some $b\in E$ over the ``boundary'' $bd(E)$, which consists  of all simplices in St$(E)$ which do not contain $E$. 
\end{definition}

\begin{theorem}
Let $X$ be a compact  topological stratified $3$--pseudomanifold. Any two triangulations of $X$ are equivalent in the sense that any triangulation of $X$ can be obtained from any other by an ambient isotopy and a finite sequence of Alexander star moves and their inverses. 
\end{theorem}

\begin{proof}
We first claim that any two triangulations of $X$ are combinatorially equivalent. In other words, give any two triangulations $k\co  X\to |K|$ and $\ell\co X\to |L|$, there exists an  isotopy $h_t\co  X \to X$ such that $h_0\co  X\to X$ is the identity and $\ell h_1k^{-1}$ is piecewise linear. In particular, there exist subdivisions $K'$ of $K$ and $L'$ of $L$ such that the induced map  $\ell h_1k^{-1}\co  K'\to L' $ is simplicial \cite[Theorem 2.14]{RS}. This claim follows from a theorem of Edward Brown \cite{Brown}:

\begin{theorem}[Brown]\label{T: Brown}
Let $K$ and $L$ be complexes of dimension $\leq 3$, and let $f\co |K|\to |L|$ be a homeomorphism. Let $\epsilon(x)$ be a positive continuous function on $K$, and let $\rho$ be a metric on $|L|$. Suppose $f$ is piecewise linear on a possibly empty subcomplex $M$ of $K$. Then there exists an isotopy $f_t\co  |K| \to |L|$ so that 
\begin{enumerate}
\item $f_0=f$; $f_1$ is piecewise linear,
\item $f_t|_{|M|}=f|_{|M|}$,
\item $\rho(f_t(x), f(x))<\epsilon(x)$ for all $x\in |K|$. 
\end{enumerate}
\end{theorem} 

We can then apply this theorem to the homeomorphism $\ell h_0k^{-1}=\ell k^{-1}\co  |K| \to |L|$ to obtain an isotopy $f_t\co  |K|\to |L|$ from $\ell k^{-1}$ to a piecewise linear map. To obtain an isotopy on $X$, define $h_t=   \ell^{-1}f_tk$. Then $h_0=\ell^{-1}f_0k=\ell^{-1}\ell k^{-1}k=\text{id}_X$ and $\ell h_1k^{-1}=\ell \ell^{-1}f_1 k k^{-1}=f_1$. 

In particular, from now on it will be convenient, with no loss of topological generality, to identify $X$ with a fixed triangulation, say $L$, and take $\ell=\text{id}\co  X\to |L|$. Then perhaps the best way to think about this theorem is to consider $k^{-1}|K|$ as giving a ``curvy'' triangulation of $X=|L|$ with respect to the piecewise linear structure induced on $X$ by $L$. Then the theorem can be interpreted as saying that there is an arbitrarily small isotopy of $X$ that straightens out the curvy triangulation given by $K$ to one that is compatible with the PL structure  given by $L$. Technically, this new ``linearized'' triangulation by $K$ is given by a new pair $(K, \bar k)$, where $\bar k\co  X\to K$ is given by $k$ precomposed with an isotopy of $X$, but it is more convenient to do without explicit use of the maps and now to imagine $K$ and $L$ as different complexes in some euclidean space which have the same polytope $|K|=|L|=X$. Piecewise linear topology then tells us that $X$ must be the polytope of a complex, say $J$, which is a subdivision of both $L$ and this straightened out $K$. 

The next step is to see that, in this context, we can transform $K$ to $L$ by a finite sequence of Alexander star moves. 

This is the content of Theorem 3.1.A of \cite{TV}:
\begin{theorem}[Turaev--Viro]\label{T: alexander}
For any polyhedron $P$ which is dimensionally homogeneous (i.e.\ it is a union of some collection of closed simplices of the same dimension) and with subpolyhedron $Q$, any two triangulations  of $P$ that coincide on $Q$ can be transformed one to another by a finite sequence of Alexander star moves and their inverses, which do not change the triangulation of $Q$.
\end{theorem}

Note that in the statement of this theorem, ``triangulation'' is being used solely in terms of triangulations already compatible with the PL structure. We have observed that any $3$--pseudomanifold is dimensionally homogeneous in Proposition \ref{P: classic def}.

Putting the above results together, we see that if we are given two, possibly piecewise linearly incompatible, triangulations of the pseudomanifold $X$, then there is first an isotopy  which makes the first triangulation piecewise linearly compatible with the second (establishing combinatorial equivalence) and then a finite sequence of Alexander star moves and their inverses which transforms the isotope of the first triangulation to the second. 
\end{proof}

Next we define Pachner's bistellar moves, which have a certain technical advantage over Alexander moves. Although these moves can be defined in arbitrary dimensions, we limit ourselves to the $3$--dimensional cases.

\begin{definition}
A \emph{bistellar move on a tetrahedron} $ABCD$  in an (abstract) simplicial complex adds a vertex $O$ at the center of $ABCD$ and replaces $ABCD$ with the four tetrahedra $OABC$, $OABD$, $OACD$ and $OBCD$. A \emph{bistellar move on a triangle} $BCD$ that is a face of exactly two tetrahedra $ABCD$ and $A'BCD$ replaces these two tetrahedra with the three tetrahedra $AA'BC$, $AA'BD$, and $AA'CD$. This move adds no new vertices but forms the new edge $AA'$. This must be viewed as a move on the abstract simplicial complex, as it may not be possible to perform this move on a complex embedded in a euclidean space. Two simplicial complexes that can be joined by a series of these two bistellar moves and their inverses are called \emph{bistellar equivalent}. 
\end{definition}

\begin{proposition}\label{P: bistellar equivalence}
Any two triangulations of a 3--pseudomanifold which are equivalent on $X(1)$ are bistellar equivalent.
\end{proposition}
\begin{proof}
By two triangulations being ``equivalent on X(1)'', we mean that the induced triangulations on $X(1)$ differ by an isotopy. This will, of course, be the case if corresponding components  of $X(1)-X(0)$ contain the same numbers of vertices. These isotopies can always be extended to $X$ as in the proof of Proposition \ref{P: realize}. So we can assume that the two triangulations actually agree on $X(1)$.  By Theorem \ref{T: alexander}, any two such triangulations are equivalent (after a further isotopy) by a sequence of Alexander star moves and inverses which do not change the triangulation of $X(1)$. This implies, in particular, that they are related by star moves (and their inverses) on simplices with interiors in $X-X(1)$. So, it suffices to show that $K$ and $\sigma_A(K)$ are bistellar equivalent, where $K$ is a triangulation of $X$ respecting the given triangulation of $X(1)$ and $\sigma_A(K)$ is the result of a star move on a simplex  $A$ with interior in $X-X(1)$. But neither star moves nor bistellar moves affect the triangulation on $\bd $St$(A)$ or exterior to St$(A)$. So consider St$(A)$ in $K$. Since  $A$ is not in $X(1)$, it follows that St$(A)\cong D^3$. But now by \cite{Cas}, $\text{St}(A)$ and $\sigma_A\text{St}(A)$ are bistellar equivalent, since this is a manifold with boundary pair. Replacing this star back into the whole complex, we see that the complexes $K$ and $\sigma_A(K)$ are bistellar equivalent.
\end{proof}

\section{Quantum invariants}\label{S: quantum}

We now show how the triangulability properties of $3$--pseudomanifolds, as studied above, can be used to obtain a generalization of  the state-sum invariants of Turaev and Viro on  $3$--manifolds to  families of invariants of $3$--pseudomanifolds, indexed by the triangulations of the natural $1$--skeleton $X(1)$. Since the natural $1$--skeleton is a compact graph, its triangulations are countable, so we obtain a countable family of invariants for each pseudomanifold. If a $3$--pseudomanifold has only point singularities, we obtain a single topological invariant.

The definition of our invariants will follow the general version given by Turaev in his monograph \cite{Tu}.  We refer the reader to this  very accessible source for the necessary definitions in the language of modular categories \cite{Tu}. In particular, we will assume a fixed unimodular category $\V$ with ground ring $K$,  braiding $c$, twist $\theta$, duality $*$, simple objects $V_i$, $i$ in the indexing set $I$, a fixed \emph{rank} $\mc D$ such that $\mc D^2=\sum_{i\in I} (\dim(V_i))^2$,    and two families of elements of $K$, $\{\dim'(i)\}_{i\in I}$ and $\{v_i'\}_{i\in I}$ such that $\dim'(0)=1$, $(\dim'(i))^2=\dim(V_i)$, $(v_i')^2=v_i$, and $v_{i^*}'=v_i'$, where $v_i$ corresponds to  $\theta\co  V_i\to V_i $ under the bijection $\Hom(V_i,V_i)\cong K$ that defines $V_i$ as a simple object.

We  note that examples of such unimodular categories exist. In fact, more generally, modular and ribbon categories by now have a rather rich history. The most well--known unimodular categories seem to be categories of certain finite dimensional representations of quantum groups. Such categories can also be defined geometrically via tangle categories, skein relations, and Temperley--Lieb algebras (in particular, see \cite{Tu}, \cite{KL}, and \cite{CFS}). We refer the reader to \cite{Tu} and the references contained therein for further exposition on this subject.

Within a strict unimodular category $\V$, Turaev constructs a certain normalized version of the 6j--symbols. In fact, he provides two alternative constructions. For the detailed construction, the reader is again referred to \cite{Tu}. We will suffice to note that the normalized $6j$--symbol $\7{i}{j}{k}{l}{m}{n}$ is an element of $H(j^*, l^*,m) \otimes_K H(i^*, j^*, k)\otimes_K H(i, m^*,n)\otimes_K H(j,l, n^*)$, where each $H(i,j,k)$ is a \emph{symmetrized} multiplicity module, i.e.\ a symmetrized version of the module $H^{ijk}=\Hom(\1, V_i\otimes V_j\otimes V_k)$, where $V_i, V_j, V_k$ are simple objects of $\V$. A symbol such as $i^*$ signifies the dual object $V_i^*$. We also note that these normalized $6j$--symbols satisfy versions of the classical identities for $6j$--symbols, such as the Biedenharn--Elliott identity. 

The $6j$--symbols can be computed explicitly in a number of cases, for example if $\V$ is the purification of the
category of representations of the quantum group $U_q(sl_2(\C))$ with
$q=-a^2$, $a$ a primitive $4r$-th root of unity in $\C$, $r\geq 2$. Formulae are provided in, e.g. \cite{Tu} and \cite{CFS}.

Now, let us return to pseudomanifolds. We define our invariant  following Turaev's for manifolds, but we must make a few generalizations. Let $X$ be a compact triangulated 3--pseudomanifold, and let $\Gamma$ be a triangulation (up to isotopy) of the natural $1$--skeleton $X(1)$ of $X$. Each edge $e$ of $X$ that does not lie in $X(1)$ possesses two normal orientations. We can think of these as the two opposite orientations of the link circle of $e$ in $X$ (or, equivalently, as Turaev does, as the two orientations of the trivial normal vector bundle to the interior of $e$ in $X$). 
If $e$ is an edge in $X(1)$, then the link of $e$ is a disjoint union of circles, and we can define a \emph{local system of orientations} for $e$ by specifying an orientation on each circle. 
If the link of $e$ consists of $n$ circles, then there are $2^n$ local systems of orientations for $e$. Let us define the \emph{set of oriented link circles of $X$} to be the set of oriented circles each of which is a component of a link of an edge of $X$. So each edge whose link has $n$ connected components makes a contribution of $2n$ elements to the set of oriented link circles (two possible orientations for each circle). To each oriented link circle $S$, there corresponds the link circle $S^*$ with the reverse orientation. 

Now, we define a \emph{coloring} $\phi$ of $X$ as an assignment from the set of oriented link circles  of $X$ to the index set $I$ of $\V$ such that if  $\phi(S)=i$, then  $\phi(S^*)=i^*$, and such that if $S_1$ and $S_2$ are both links of the same edge then $\phi(S_1)$ equals either $\phi(S_2)$ or $\phi(S_2)^*$. Let col$(X)$ denote the set of coloring of $X$. If $X$ is a manifold, then this set of colorings is equivalent to the set of normally oriented edge colorings of \cite{Tu}. 

For each edge $e$ of $X$ with link $S$,  let $\dim_{\phi}(e)=\dim(\phi(S))=\dim(\phi(S^*))$. This is well--defined by our definition of $\phi$. 

Suppose that $T$ is a $3$--simplex of $X$ and that we are given a coloring $\phi$ of $X$. Let us arbitrarily label the vertices of $T$ by $A$, $B$, $C$, and $D$.  Then for each edge $e$ of $T$, there is a unique (non--oriented) link circle $S_{T,e}$ of $e$ that intersects $T$ for a small enough choice of distinguished neighborhood of a point in  $e$. Let $i$ be the color of $S_{T,AB}$ under $\phi$ with orientation determined to agree with the direction from $\vec{AD}$ to $\vec{AC}$ within $T$. Similarly, let $j,k,l,m,n$ be the colors associated by $\phi$ to the links $S_{T,BC}$, 
$S_{T,AC}$, $S_{T,CD}$, $S_{T,AD}$, and $S_{T,BD}$ with orientations determined by the pairs $(\vec{BD},\vec{BA})$, $(\vec{AB},\vec{AD})$, $(\vec{CB},\vec{CA})$, $(\vec{AC},\vec{AB})$, and $(\vec{BA},\vec{BC})$.
 We then define $|T^{\phi}|=\7{i}{j}{k}{l}{m}{n}$.

Now, given a coloring $\phi\in$col$(X)$, let $$|X|^{\Gamma}_{\phi}=\mc D^{-2a}(\prod_e \dim_{\phi}(e))\text{cntr}(\otimes_T |T^{\phi}|),$$
where $a$ is the number of vertices in the triangulation of $X$, $e$ runs over all (non-oriented) edges of $X$, $T$ runs over all $3$--simplices of $X$, and $\text{cntr}$ indicates tensor contraction. This contraction is well--defined since wherever two tetrahedra share a face in the triangulation, the corresponding $|T|$'s will have tensor factors of the respective forms $H(i,j,k)$ and $H(i^*,j^*,k^*)$ (note that for two such neighboring tetrahedra $T$ and $T'$, the link determined for any common edge $e$ is the same for each tetrahedron, i.e. $S_{T,e}=S_{T',e}$). These modules have a non--degenerate bilinear pairing to $K$ induced by the trace, and we use this pairing to contract. Since all pairs of faces contract this way, $|X|_{\phi}$ is an element of $K$. 

Finally, we define the state sum $$|X|^{\Gamma}=|X|^{\Gamma}_{\V}=\sum_{\phi\in \text{col}(X)}|X|^{\Gamma}_{\phi}\in K.$$

\begin{theorem}
Let $X$ be a compact $3$--pseudomanifold and $\V$ a strict unimodular category. The state sum $|X|_{\V}^{\Gamma}$ is independent of the triangulation of $X$ modulo the given isotopy class of the fixed triangulation $\Gamma$ of $X(1)$. In particular, for fixed $\V$, we obtain a family of topological invariants of $X$ indexed by the (countable) triangulations up to isotopy  of the natural $1$--skeleton $X(1)$.  
\end{theorem}
\begin{proof}

By Proposition \ref{P: bistellar equivalence}, any two triangulations of $X$ that agree on $X(1)$ are bistellar equivalent, i.e.\ we can pass from one to the other by a series of isotopies and bistellar moves. Furthermore, none of these bistellar moves add or remove any edges in $X(1)$. Clearly, $|X|^{\Gamma}$ is unchanged by isotopy, so we must see that it is unchanged by bistellar moves. But we now refer the reader to Turaev's proof in \cite{Tu} that the corresponding state sum is an invariant  for (closed) manifolds. This is proven  by  applying locally identities on the normalized $6j$--symbols that generalize more classical identities on $6j$--symbols including the Biedenharn--Elliott identity. Since the edges in $X(1)$ with their more complicated colorings remain fixed under the bistellar moves, Turaev's proof applies directly. 
\end{proof}

If $X$ has only point singularities, we obtain a topological invariant:

\begin{corollary}
Let $X$ be a compact $3$--pseudomanifold with only point singularities. Then $|X|_{\V}$ is a topological invariant of $X$.
\end{corollary}

This is essentially the invariant of Barrett and Westbury \cite{BW}, though they work in a slightly different category system.

This corollary also allows us to obtain an invariant of manifolds with boundary by coning off the boundary pieces. If a manifold has multiple boundary pieces, there will be multiple ways to cone, as noted in the statement of the next corollary. However, as we will see in the following corollary, the various ways of coning give fundamentally the same result. 

\begin{corollary}\label{C: boundary}
Let $M$ be a $3$--manifold with $n$ boundary components. Let $p(n)$ denote the number of ways to partition $n$ objects into groups. Then there exist $p(n)$ topological invariants $|M|^i_{\V}$, where $i$ denotes the $i$th partition (in some given fixed ordering). These are defined by $|M|^i_{\V}=|X^i|_{\V}$, where $X^i$ is the $3$--pseudomanifold obtained from $M$ by adding one cone on each group of boundaries as determined by the partition.  (Some of the invariants may be redundant, dependent upon the symmetries of the manifold.) 
\end{corollary}

\begin{corollary}
Suppose that the compact $3$--pseudomanifold $Y$ is obtained from the
compact $3$--pseudomanifold $X$ by a quotient map $p$ that is
injective on the complement of a finite set $S$ of $n$ points  and
takes this set of $n$ points to a set of $m<n$ points in $Y$. Let
$\Gamma$ be a triangulation of $X(1)$ such that each point  $x\in
S\cap X(1)$ with $\#\{p^{-1}(p(x))\}>1$ is a vertex of $\Gamma$.
Then for any $\V$,
$|Y|^{\Gamma}_{\V}=\mc{D}^{2(n-m)}|X|^{\Gamma}_{\V}$.
\end{corollary}
\begin{proof}
We can find triangulations of $X$ and $Y$ such that the triangulation
of $Y$ is formed from the triangulation of $X$ by identifying certain
vertices. In fact, our assumption on $\Gamma$ ensures that all points
of $S$ in $X(1)$ are vertices already, and any other point in $S$ can
be made a vertex without disturbing $\Gamma$ by employing a
generalized stellar subdivision. Then we need only identify vertices
as prescribed by $p$ to obtain a triangulation of $Y$. The formula
now follows from the definitions, since these triangulations of $X$
and $Y$ give equivalent combinatorial data for the computation of 
$|\cdot|^{\Gamma}_{\V}$, except for the number of vertices, which
comes in only as the  multiplicative factor $\mc D^{-2a}$. 
\end{proof}

So we see that the various invariants of Corollary \ref{C: boundary} really only differ by  a well--controlled factor of $\mc{D}$. In fact, each of these invariants must be a contraction of the vector obtained from the manifold with boundary by considering the topological quantum field theory version of the Turaev--Viro theory.

Finally, we indicate how our state sum invariants of pseudomanifolds can be used to yield invariants of tame knots and links. This is done by assigning to knots and links certain pseudomanifolds. 

The simplest such construction is the following: Given a knot or link in $S^3$, we can consider the complement of an open regular or tubular neighborhood of the knot or link. This complement is a $3$--manifolds whose boundary consists of tori. We can then cone off collections of these boundary tori to obtain pseudomanifolds as in Corollary \ref{C: boundary}.  By the same method, we can obtain invariants of links in any $3$--manifold (if the manifold is not orientable, we may have to replace cones on tori with cones on Klein bottles). We can even apply this process to links in a $3$--pseudomanifold $X$ since, for a PL link, its regular neighborhood will remain well--defined, thinking of $X$ as a $PL$ space. Again, we can remove an open regular neighborhood and then cone off the collection of free $2$--simplices. 
Most generally, this method also provides invariants of pairs $(X,A)$, where $X$ is any dimensionally homogeneous $3$--complex and $A$ is a subcomplex. 

We can also create from knots and links pseudomanifolds with $1$--dimensional singularities: If $\bar N=S^1\times D^2$ represents the closed regular neighborhood of a knot in an orientable $3$--manifold $M$, $N$ its interior, and $\bd \bar N=S^1\times S^1$ its boundary, define $f_k\co S^1\times S^1\to S^1$ by $(\theta, \phi)\to k\theta$. Let $C_k$ be the mapping cylinder of $f_k$. Then $(M-N)\cup_{\bd \bar N} C_k$ is a pseudomanifold with a $1$--dimensional singular set. For links, we can either perform this construction for each component individually or, for  $n$--component sublinks,  adjoin the mapping cylinders of the maps $f_{(k_1, \ldots, k_n)}\co  \coprod_{i=1}^n S^1\times S^1\to S^1$, where $f_{(k_1, \ldots, k_n)}$ restricted to the boundary of the neighborhood of the $i$th component of the link is $f_{k_i}$. If $M$ is not orientable, we can still perform this procedure after replacing some of the tori with Klein bottles; the maps $f_k$ can still be defined since they factor through the projection to the base of the circle bundle.

Finally, we can, of course, employ the constructions of the previous two paragraphs in combinations. Hence, for links, we obtain several families of associated pseudomanifolds and thus several families of invariants.

\providecommand{\JFM}[1]{\relax\ifhmode\unskip\space\fi
  \href{http://www.emis.de/cgi-bin/JFM-item?#1}{JFM~#1}}
\providecommand{\ZBL}[1]{\relax\ifhmode\unskip\space\fi
  \href{http://www.emis.de/cgi-bin/zmen/ZMATH/en/quick.html?type=html&an=#1}{Zbl~#1}}

\def\Dummy{\hbox to 1in{$\phantom{\rm My}$\hss}}

\Addresses\recd


\begin{thebibliography}

\bibitem{Alex}
\textbf{J\,W Alexander}, \href{http://links.jstor.org/sici?sici=0003-486X(193005)2:31:2%3C292:TCTOC%3E2.0.CO%3B2-T}{\Dummy}\kern -1in\emph{The combinatorial theory of complexes}, Ann. of
  Math. 31 (1930) 292--320, \MR{1 502 943}, \JFM{56.0497.02}

\bibitem{BW} \textbf{John~W Barrett}, \textbf{Bruce~W Westbury},
\href{http://www.ams.org/jourcgi/jour-getitem?pii=S0002994796016601}{\Dummy}\kern -1in\emph{Invariants of
piecewise-linear 3-manifolds}, Trans. Amer. Math. Soc. 348 (1996)
3997--4022, \MR{97f:57017}, \ZBL{0865.57013}

\bibitem{Bo}
\textbf{A Borel~et al}, \emph{Intersection Cohomology}, volume~50 of
  \emph{Progress in Mathematics}, Birkhauser, Boston (1984),
  \MR{88d:32024}, \ZBL{0553.14002}

\bibitem{Brown}
\textbf{Edward~M Brown}, \href{http://links.jstor.org/sici?sici=0002-9947%28196910%29144%3C173%3ATHF3%3E2.0.CO%3B2-%23}{\Dummy}\kern -1in\emph{The {H}auptvermutung for 3--complexes}, Trans.
  Amer. Math. Soc. 144 (1969) 173--196, \MR{40:4956},
  \ZBL{0198.28103}

\bibitem{CFS} \textbf{J\,Scott Carter}, \textbf{Daniel~E Flath},
\textbf{Masahico Saito}, \emph{The Classical and Quantum 6j--symbols},
Mathematical Notes, vol.~43, Princeton University Press, Princeton, NJ
(1995) \MR{97g:17008}, \ZBL{0851.17001}

\bibitem{Cas}
\textbf{Maria~Rita Casali}, \emph{A note about bistellar operations on
  {PL}-manifolds with boundary}, Geom. Dedicata 56 (1995) 257--262,
\MR{96h:57021}, \ZBL{0833.57008}

\bibitem{Co}
\textbf{Marshall~M Cohen}, \href{http://links.jstor.org/sici?sici=0002-9947%28196902%29136%3C189%3AAGTORR%3E2.0.CO%3B2-E}{\Dummy}\kern-1in\emph{A general theory of relative regular
  neighborhoods}, Trans. Amer. Math. Soc. 136 (1969) 189--229, \MR{40:2052},
  \ZBL{0182.57602}

\bibitem{GM1}
\textbf{Mark Goresky}, \textbf{Robert MacPherson}, \href{http://dx.doi.org/10.1016/0040-9383(80)90003-8}{\Dummy}\kern -1in\emph{Intersection Homology
  Theory}, Topology 19 (1980) 135--162, \MR{82b:57010}, \ZBL{0448.55004}

\bibitem{GM2}
\textbf{Mark Goresky}, \textbf{Robert MacPherson}, \emph{Intersection Homology
  {II}}, Invent. Math. 72 (1983) 77--129, \MR{84i:57012}, \ZBL{0529.55007}

\bibitem{Ha}
\textbf{Allen Hatcher}, \href{http://www.math.cornell.edu/~hatcher/AT/AT.pdf}{\Dummy}\kern -1in\emph{Algebraic Topology}, Cambridge University Press,
  Cambridge (2002)

\bibitem{HUD}
\textbf{J\,F\,P Hudson}, \emph{Piecewise Linear Topology}, W\,A Benjamin, Inc. New
  York (1969), \MR{40:2094},
  \ZBL{0189.54507}

\bibitem{KL} \textbf{Louis~H Kauffman}, \textbf{S\'ostenes Lins},
\emph{Temperley-{L}ieb Recoupling Theory and Invariants of
3-Manifolds}, Annals of Math. Studies 134, Princeton University
Press, Princeton, NJ (1994), \MR{95c:57027},
  \ZBL{0821.57003}

\bibitem{Mo52}
\textbf{Edwin~E Moise}, \href{http://links.jstor.org/sici?sici=0003-486X(195207)2:56:1%3C96:ASI3VT%3E2.0.CO%3B2-2}{\Dummy}\kern -1in\emph{Affine structures in 3-manifolds {V}: The
  triangulation theorem and {H}auptvermutung}, Ann. of Math. 56 (1952) 96--114

\bibitem{MK3}
\textbf{J Munkres}, \emph{The triangulation of locally triangulable spaces},
  Acta Math 97 (1957) 67--93, \MR{19,437b}, \ZBL{0080.16801}

\bibitem{MK2}
\textbf{James~R Munkres}, \emph{Topology: A First Course}, Prentice-Hall, Inc.,
  Englewood Cliffs, NJ (1975), \MR{57:4063}, \ZBL{0306.54001}

\bibitem{MK}
\textbf{James~R Munkres}, \emph{Elements of Algebraic Topology},
  Addison-Wesley, Reading, MA (1984), \MR{85m:55001}, \ZBL{0673.55001}

\bibitem{Pa}
\textbf{U Pachner}, \emph{{B}istellare {{\"{A}}quivalenz} kombinatorischer
  {M}annigfaltigkeiten}, Arch. Math 30 (1978) 89--98, \MR{58:7641}, \ZBL{0375.57007}

\bibitem{RS} \textbf{C\,P Rourke}, \textbf{B\,J Sanderson},
\emph{Introduction to Piecewise-Linear Topology}, Springer Study
Edition, Springer-Verlag, Berlin-Heidelberg-New York (1982) \MR{83g:57009},
\ZBL{0477.57003}

\bibitem{Tu} \textbf{V\,G Turaev}, \emph{Quantum Invariants of Knots
and 3-Manifolds}, de Gruyter Studies in Mathematics 18, Walter de
Gruyter, Berlin-New York (1994) \MR{95k:57014}, \ZBL{0812.57003}

\bibitem{TV}
\textbf{V\,G Turaev}, \textbf{O\,Y Viro}, \emph{State sum invariants of 3-manifolds
  and quantum 6j-symbols}, Topology 31 (1992) 865--902, \MR{94d:57044}, \ZBL{0779.57009}

\end{thebibliography}
\end{document}